\documentclass[english]{amsart}
\usepackage[T1]{fontenc}
\usepackage[latin1]{inputenc}
\usepackage{babel}
\usepackage{graphics}

\makeatletter

 \theoremstyle{plain}    
 \newtheorem{thm}{Theorem}[section]
 \numberwithin{equation}{section} 
 \numberwithin{figure}{section} 
 \theoremstyle{plain}    
 \newtheorem{cor}[thm]{Corollary} 
 \theoremstyle{plain}    
 \newtheorem{lem}[thm]{Lemma} 
 \theoremstyle{plain}    
 \newtheorem{prop}[thm]{Proposition} 
 \theoremstyle{definition}
 \newtheorem{defn}[thm]{Definition}
 \theoremstyle{remark}
 \newtheorem{rem}[thm]{Remark}
 \theoremstyle{remark}
 \newtheorem*{rem*}{Remark}
 \theoremstyle{remark}    
 \newtheorem{notation}[thm]{Notation} 

\makeatother
\begin{document}

\title{An involution on the quantum cohomology ring of the grassmannian}

\author{Harald Hengelbrock}



\curraddr{Fakultät für Mathematik, Ruhr-Universität Bochum, 44803 Bochum}

\email{harald.hengelbrock@ruhr-uni-bochum.de}

\date{February 14, 2002}

\begin{abstract}
For a Fano manifold M, complex conjugation defines a real involution
on the quantum cohomology ring \( QH^{*}(M,\Bbb {C}) \). For the
Grassmannian we identify this involution with an explicit transformation
on Schubert classes defined over the integers. It is a composition
of a power of the cyclic action recently defined by Agnihotri and
Woodward with Poincaré duality. The existence of the involution has
been observed independently and from a different point of view by
Postnikov (math.CO/0205165).
\end{abstract}
\maketitle

\section{Introduction}

The quantum cohomology ring of the Grassmann varieties has been the
subject of many recent studies, since it is related to different branches
of mathematics, and thus displays a rich mathematical structure. One
of its features is its semisimplicity, that is, as an abstract algebra
over the complex numbers it is isomorphic to the coordinate ring of
a finite set of reduced points. 

Now any permutation of the closed points of \( \textrm{Spec}\, R \)
induces an automorphism of the quantum ring. It turns out that the
natural permutation given by complex conjugation of the points of
\( \textrm{Spec}\, R \) has a meaning in terms of Schubert classes.
Indeed, it even defines an involution of the quantum ring over the
integers. This involution is seen to be strongly related to Poincaré
duality of the classical cohomology ring. The surprising fact is that,
combined with an action recently defined by Agnihotri and Woodward
in \cite{AW}, Poincaré duality is a ring automorphism.

In Theorem \ref{th1} we prove that a given explicit transformation
of Young diagrams induces an involution on the quantum ring. The second
part is devoted to the proof that this involution is nothing but complex
conjugation on the function level.

\section{The involution}

For basic facts about the small quantum cohomology ring of Grassmannians
we refer to \cite{BIF} and the references therein.

\begin{notation}
Let \( G=G(k,n) \) denote the Grassmannian of \( k \)-planes in
\( \Bbb {C}^{n}. \) Set \( l:=n-k. \) For a Schubert class \( S \)
let \( \hat{S} \) denote the Poincaré dual of \( S \). This operation
extends additively to general classes. \( H^{*}(G,\Bbb {Z}) \) is
the classical cohomology ring, while \( QH^{*}(G,\Bbb {Z}) \) denotes
the quantum cohomology ring of \( G \) with coefficients in \( \Bbb {Z} \),
where the constant \( q \) of the standard definition is set to \( 1 \).
Products in \( QH^{*}(G,\Bbb {Z}) \) will be denoted by a `\( * \)',
while `\( \cup  \)' is used for the classical product. Three point
invariants will be written in the form \( \langle A,B,C\rangle  \)
for Schubert classes \( A,B,C \). 

A Schubert class \( S_{\lambda } \) is often represented by its Young
diagram \( \lambda =(\lambda _{1},\ldots ,\lambda _{l}) \). 
\end{notation}
\begin{defn}
Define an involution \( S\rightarrow \overline{S} \) on the \( H^{*}(G,\Bbb {Z}) \)
by the following rule:\\
Let \( \lambda =(\lambda _{1},\ldots ,\lambda _{l}) \) be a Young
diagram corresponding to a Schubert class \( S \). Let \( d_{\lambda } \)
be the largest \( i \) with \( \lambda _{i}\geq i. \) Then \( \overline{S} \)
is the Schubert class corresponding to the Young diagram \( \mu =(\mu _{1},\ldots ,\,  \)\( \mu _{l}) \)
defined by \[
\mu _{i}=\left\{ \begin{array}{ccc}
d_{\lambda }+k-\lambda _{d_{\lambda }-i+1} & : & i\leq d_{\lambda }\\
d_{\lambda }-\lambda _{l-i+d_{\lambda }+1} & : & i>d_{\lambda }.
\end{array}\right. \]

Graphically, \( d_{\lambda } \) is the width of the biggest square
inside \( \lambda  \), and \( \mu  \) is obtained by dualizing the
sub-diagrams in the upper left and lower right rectangles in the subdivision
given by drawing lines on the outer border of the square. 
\end{defn}
\vspace{0.3cm}
{\centering \includegraphics{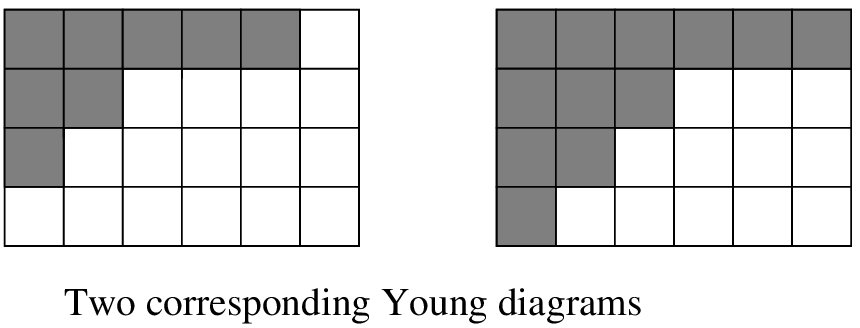} \par}
\vspace{0.3cm}

Again, this operation extends linearly to general classes. 

As checked in the following Lemma, this definition is equal to the
operator \( S\rightarrow \widehat{c^{k}S} \), where \( c^{k} \)
is defined as in Chapter 7 of \cite{AW}:

The operator \( c \) is given by quantum multiplication with the
Schubert class corresponding to the Young diagram \( C=(1,1,\ldots ,1) \),
and \( c^{k} \) correspondingly by multiplying with the \( k \)-th
power of \( C \). Powers of \( c \) give a \( \Bbb {Z}/n \) action
on \( QH^{*}(G,\Bbb {Z}) \) (\cite{AW}, Proposition 7.2.). 

\begin{rem*}
\( C^{k} \) is represented by the diagram \( (k,\ldots ,k), \) which
is the class of a point on \( G(k,n). \) 
\end{rem*}
\begin{lem}
The operator \( S\mapsto \overline{S} \) equals \( S\mapsto \widehat{c^{k}S.} \) 
\end{lem}
\begin{proof}
For a Schubert class \( S \) given by \( (\lambda _{1},\ldots ,\lambda _{l}) \),
set \( \lambda ^{i}=k+i-\lambda _{i}. \) In the notation used in
\cite{AW}, the Schubert class is then expressed by \( (\lambda ^{1},\ldots ,\lambda ^{l}) \).
In this notation, \( c^{k}S \) is represented by \[
\eta =(\lambda ^{d_{\lambda }+1}-k,\ldots ,\lambda ^{l}-k,\lambda ^{1}-k+n,\ldots ,\lambda ^{d_{\lambda }}-k+n)\]
 (\cite{AW}, beginning of Chapter 7, the number \( m \) is exactly
\( d_{\lambda }, \) and \cite{AW}, Proposition 7.2.). As a Young
diagram, \( c^{k}S \) is represented by \( (\mu _{1},\ldots ,\mu _{l}) \)
with\[
\mu _{i}=k+i-\eta _{i}=\left\{ \begin{array}{ccc}
k-d_{\lambda }+\lambda _{d_{\lambda }+i} & : & i\leq l-d_{\lambda }\\
-d_{\lambda }+\lambda _{d_{\lambda }-l+i} & : & i>l-d_{\lambda }.
\end{array}\right. \]
 But this is just the Young diagram of the Poincaré dual of \( \overline{S}. \)
\end{proof}
\begin{thm}
\label{th1}The involution operator defined above respects the quantum
product: for \( A,B\in QH^{*}(G,\Bbb {Z}) \) the formula \( \overline{A*B}=\overline{A}*\overline{B} \)
holds.
\end{thm}
\begin{proof}
Due to associativity and linearity it is enough to show this equation
for a Schubert class \( A \) and \( B=\sigma _{r}, \) where \( \sigma _{r} \)
corresponds to the diagram \( (r,0,\ldots ,0), \) since the \( \sigma _{i} \)
generate the quantum ring. Using the definition of the quantum product,
write \[
\overline{A*B}=\sum \langle A,B,\widehat{T_{i}}\rangle \overline{T_{i}}\]
 where the sum runs over the Schubert basis, and\[
\overline{A}*\overline{B}=\sum \langle \overline{A},\overline{B},\widehat{T_{i}}\rangle T_{i}=\sum \langle \overline{A},\overline{B},\widehat{\overline{T_{i}}}\rangle \overline{T_{i}}.\]
Comparing the coefficients we find that the claim is equivalent to
the following equation of invariants:\[
\langle A,\sigma _{r},\hat{S}\rangle =\langle \overline{A},\overline{\sigma _{r}},\widehat{\overline{S}}\rangle \]
 for all \( A,S,r. \) But\[
\langle \overline{A},\overline{\sigma _{r}},\widehat{\overline{S}}\rangle =(\overline{A}*\widehat{\overline{S}})\cup \overline{\sigma _{r},}\]
and \( \overline{A}*\widehat{\overline{S}}=\widehat{c^{k}A}*c^{k}S=c^{n-k}\hat{A}*c^{k}S=\hat{A}*c^{n}S=\hat{A}*S. \)
Here the third equality comes from the commutativity of the quantum
product and the second equality comes from the fact that \( \overline{\overline{A}}=A, \)
which implies \begin{equation}
\label{eq1}
\widehat{c^{k}A}=c^{n-k}\widehat{A}.
\end{equation}

We see that \( \langle \overline{A},\overline{\sigma _{r}},\widehat{\overline{S}}\rangle =\langle \hat{A},S,\overline{\sigma _{r}}\rangle , \)
so it remains to show that\begin{equation}
\label{for1}
\langle A,S,\sigma _{r}\rangle =\langle \hat{A},\hat{S},\overline{\sigma _{r}}\rangle .
\end{equation}
 Let A and S be given by Young diagrams \( (a_{1,}\ldots ,a_{l}) \)
and \( (s_{1},\ldots ,s_{l}). \) The (quantum) Pieri formula states
that \( \langle A,S,\sigma _{r}\rangle  \) is zero except for two
cases (\cite{BIF}, Example 2 and \cite{Bu}, Section 4), where it
takes value \( 1 \):

\begin{equation}
\label{eq2}
\textrm{deg }A+\textrm{deg }S+r=kl\textrm{ with }\left\{ \begin{array}{ccc}
a_{i}+s_{j}\geq k & : & i+j=l\\
a_{i}+s_{j}\leq k & : & i+j=l+1
\end{array}\right. 
\end{equation}
 \begin{equation}
\label{eq3}
\textrm{deg}\; A+\textrm{deg}\; S+r=kl+n\textrm{ with }\left\{ \begin{array}{ccc}
a_{i}+s_{j}\geq k+1 & : & i+j=l+1\\
a_{i}+s_{j}\leq k+1 & : & i+j=l+2
\end{array}\right. 
\end{equation}

If these degree conditions do not hold, the same will be true for
\( \widehat{A} \) , \( \widehat{B} \) and \( \overline{\sigma _{r}} \)
. We can therefore assume \( \textrm{deg }A+\textrm{deg }S+r=kl+an,\, a\in \{0,1\}. \) 

We can write \( \overline{\sigma _{r}} \), given by the diagram \( (k-r,1,\ldots ,1), \)
as \( c^{1}\sigma _{k-r}. \) The definition of \( c \) and commutativity
yield \[
\langle \hat{A},\hat{S},\overline{\sigma _{r}}\rangle =\langle \hat{A},\hat{S},c^{1}\sigma _{k-r}\rangle =\langle \hat{A},c^{1}\hat{S},\sigma _{k-r}\rangle .\]

In case (\ref{eq2}) it follows that \( \textrm{deg}\; \hat{A}+\textrm{deg}\; \hat{S}+\textrm{deg}\; \overline{\sigma _{r}}=(kl-\textrm{deg}\, A)+(kl-\textrm{deg}\, S)+(n-r)=2kl+n-kl=kl+n, \)
while the conditions on \( A \) and \( S \) become \begin{equation}
\label{eq4}
\left\{ \begin{array}{ccc}
a_{i}'+s_{j}'\geq k & : & i+j=l+1\\
a_{i}'+s_{j}'\leq k & : & i+j=l+2,
\end{array}\right. 
\end{equation}
 where the \( a_{i}' \) and \( s_{i}' \) give the diagrams for \( \hat{A} \)
and \( \hat{S} \). 

Let us now look at \( s_{1}'. \) If \( s_{1}'<k \), the dual Pieri
formula (\cite{BIF}, Proposition 4.2) yields \( (c^{1}\widehat{S})_{i}=s'_{i}+1, \)
where \( (c_{1}\widehat{S})_{i} \) means the \( i-\textrm{th} \)
component of the Young diagram corresponding to \( c^{1}\widehat{S} \).
Our conditions (\ref{eq4}) become exactly the conditions of the Pieri
formula for the invariant \( \langle \hat{A},c^{1}\hat{S},\sigma _{k-r}\rangle  \)
to be nonzero. This implies equation (\ref{for1}).

If \( s_{1}'=k, \) the quantum dual Pieri formula yields \( (c^{1}\widehat{S})_{i}'=s_{i+1}' \).
The conditions (\ref{eq4}) then result in the classical Pieri conditions
for \( \langle \hat{A},c^{1}\hat{S},\sigma _{k-r}\rangle  \), so
that equality (\ref{for1}) holds.

The case (\ref{eq3}) implies \( \textrm{deg}\; \hat{A}+\textrm{deg}\; \hat{S}+\textrm{deg}\; \overline{\sigma _{r}}=kl, \)
and the conditions\[
\left\{ \begin{array}{ccc}
a_{i}'+s_{j}'\geq k-1 & : & i+j=l\\
a_{i}'+s_{j}'\leq k-1 & : & i+j=l+1.
\end{array}\right. \]
 Write again \( \langle \hat{A},\hat{S},\overline{\sigma _{r}}\rangle =\langle \hat{A},c^{1}\hat{S},\sigma _{k-r}\rangle  \)
and note that \( s'_{1}<k \) by the first condition on the \( a_{i} \)
and \( s_{i} \) in (\ref{eq3}) above. As above, this implies \[
(c^{1}\widehat{S})_{i}=s'_{i}+1,\]
 so that, as in the first case, the statement follows by the Pieri
formula.
\end{proof}
\begin{cor}
\label{cor1}For any classes \( A,B,C\in QH^{*}(G,\Bbb {Z}) \)  the
following equalities hold:\[
\widehat{(A*C)}=\widehat{A}*\overline{C}\textrm{ }\]
\[
\langle A,C,B\rangle =\langle \widehat{A},\widehat{C},\overline{B}\rangle .\]

\end{cor}
\begin{proof}
Consider \( \widehat{A*C}=\widehat{A*C}*c^{n-k}*c^{k}=\overline{A*C}*c^{k}=c^{k}\overline{A}*\overline{C}=\widehat{A}*\overline{C} \),
where we use \ref{eq1} and the Theorem. In terms of GW-invariants,
this results in \[
\sum \langle A,C,T_{i}\rangle T_{i}=\sum \langle \widehat{A},\overline{C},\widehat{T_{i}}\rangle T_{i},\]
 where \( T_{i} \) ranges over the Schubert classes. We obtain the
second equality by comparing the invariants and expressing \( B \)
in the Schubert basis.
\end{proof}
We now consider the quantum cohomology ring over the complex numbers.
Set \( R:=QH^{*}(G,\Bbb {C}) \). The space \( \textrm{Spec}\: R \)
consists of \( \textrm{dim}\: R \) reduced points, that is, \( R \)
is a semisimple algebra over the complex numbers (\cite{ST}, Theorem
4.6). The points of \( \textrm{Spec}\, R \) will be considered as
embedded into complex affine space. Using the isomorphism of the local
rings with the complex numbers, we may regard the cohomology classes
as complex functions on \( \textrm{Spec}\, R \). We wish to show
that, on this level, our involution is given by complex conjugation. 

\begin{defn}
For any cohomology class \( C\in QH^{*}(G,\Bbb {Z}) \) define \( \phi _{C} \)
to be the endomorphism of the finite vector space \( QH^{*}(G,\Bbb {Z}) \)
induced by quantum multiplication with \( C. \)
\end{defn}
\begin{prop}
\label{sympos}For any cohomology class \( C\in QH^{*}(G,\Bbb {Z}) \)
, the endomorphism \( \phi _{C*\overline{C}} \) is represented by
a semipositve definite symmetric matrix with respect to the Schubert
basis. The cohomology class \( C*\overline{C} \) defines a real-valued
positive function on \( \textrm{Spec}\, R. \) 
\end{prop}
\begin{proof}
Symmetry of the map \( \phi  \) is equivalent to the following statement:
for any two Schubert classes \( A \) and \( B \) we have \[
A*C*\overline{C}\cup \widehat{B}=B*C*\overline{C}\cup \widehat{A},\]
 that is, the coefficient of \( B \) in the product \( A*C*\overline{C} \)
is the same as the coefficient of \( A \) in the product \( B*C*\overline{C}. \)
To prove this, consider\[
A*C*\overline{C}\cup \widehat{B}=\sum \langle A,C,\widehat{T_{i}}\rangle \langle T_{i},\overline{C},\widehat{B}\rangle \]
 \[
=\sum \langle \widehat{A},\overline{C},T_{i}\rangle \langle \widehat{T_{i}},C,B\rangle =B*C*\overline{C}\cup \widehat{A},\]
 where the sums run over the Schubert basis. The second and third
equalities are derived from the definition of the quantum product,
while the second uses Corollary \ref{cor1}. 

For semipositivity, given an arbitrary cohomology class K, we show
that\[
K*B*\overline{B}\cup \widehat{K}\geq 0.\]
 We write again \[
K*B*\overline{B}\cup \widehat{K}=\sum \langle K,B,\widehat{T_{i}}\rangle \langle T_{i},\overline{B},\widehat{K}\rangle =\sum \langle K,B,\widehat{T_{i}}\rangle \langle K,B,\widehat{T_{i}}\rangle \geq 0\]
 where the second equality uses Corollary \ref{cor1}.

The last statement now follows, because a complex value of \( C*\overline{C} \)
on some point \( p \) of \( \textrm{Spec}\, R \) would induce a
complex eigenvalue of \( \phi  \) with eigenvector the function taking
value 1 on \( p \) and \( 0 \) on the other points. 
\end{proof}
\begin{lem}
\label{lemns}Consider \( C \) as a complex function on \( \textrm{Spec}\, R \).
\( C \) vanishes on some point \( p \) if and only if \( \overline{C} \)
does. 
\end{lem}
\begin{proof}
Assume \( C(p)\neq 0 \) but \( \overline{C}(p)=0. \) Consider \( A:=C-\overline{C} \)
. The value of \( A*\overline{A} \) on \( p \) is \( C(p)^{2}. \)
But, since \( A \) is anti-invariant, \( A*\overline{A}=-A^{2}. \)
This implies \( C(p)^{2}=-C(p)^{2}, \) a contradiction. 
\end{proof}
\begin{lem}
\label{prop1}Let \( C \) be any cohomology class over \( Z \) ,
\( C=a+ib \) with real functions \( a,b \). Then there exists a
real positive function \( r \) with \[
\overline{C}=r(a-ib),\]
 and \( C*\overline{C} \) is a real function on \( \textrm{Spec}\, R. \) 
\end{lem}
\begin{proof}
Consider a point \( p\in \textrm{Spec}\, R \) . If \( C(p) \) is
zero, so is \( \overline{C}(p) \) by the preceeding Lemma. Otherwise,
\( C*\overline{C}(p) \) is real and positive by Proposition \ref{sympos}. 
\end{proof}
\begin{rem}
\label{rem1}The function \( r \) is nonzero wherever \( C \) is
nonzero, since \( \overline{\overline{C}}=C. \) 
\end{rem}
\begin{prop}
\label{Prop3}The function r of Lemma \ref{prop1} can be taken to
be identically 1 on \( \textrm{Spec}\, R. \) 
\end{prop}
\begin{proof}
Take a Schubert class \( S=a+ib,\, \overline{S}=r(a-ib) \) as in
Lemma \ref{prop1}. We first treat the case \( b(p)=0 \) for \( p\in \textrm{Spec}\, R \).
It is enough to show that \( B:=S-\overline{S}=a-ra=0 \) on \( p. \)
But if \( B(p)\neq 0, \) the real function \( B*\overline{B}=-B^{2} \)
would take a negative value, contradicting Proposition \ref{sympos}. 

Now set \( A:=S+\overline{S}=c+id, \) with \( c,d \) real. Then,
again by Lemma \ref{prop1}, \( \overline{A}=r'(c-id), \) with a
real function \( r'. \) First, note that \( c \) must be nonzero
on all points \( p \) where \( d \) is, since otherwise \( A*\overline{A}=A^{2}=-d^{2} \)
at \( p \), contradicting Proposition \ref{sympos}. But \( A=\overline{A}, \)
so that \( c+id=r'(c-id). \) It follows that \( d \) is always zero,
and that \( r' \) is identically 1 where \( A \) is nonzero. Thus,
for points \( p\in \textrm{Spec}\, R \) with \( S(p)\neq 0, \) from
the identity \( d=b-rb=0 \) with \( b(p)\neq 0 \) we deduce that
\( r(p)=1. \) 
\end{proof}
\begin{thm}
The involution \( S\rightarrow \overline{S} \) is, at the function
level, equal to complex conjugation. 
\end{thm}
\begin{proof}
Immediate from Lemma \ref{prop1} and Proposition \ref{Prop3}.
\end{proof}

\end{document}